\documentclass{article}
\pdfoutput=1
\usepackage[utf8]{inputenc}
\usepackage[english]{babel}
\usepackage{amsmath, amsthm}
\usepackage[shortlabels]{enumitem}
\usepackage{tikz}
\usetikzlibrary{matrix}
\usetikzlibrary{cd}
\usepackage{amsfonts, amssymb}
\usepackage{hyperref}
\usepackage[left=4cm,right=4cm,top=2cm,bottom=2cm]{geometry}
\usepackage{csquotes}
\usepackage[style=alphabetic,backend=biber]{biblatex}

\newcommand{\C}{\mathbb{C}}
\newcommand{\Gl}{\mathrm{GL}}
\newcommand{\Sl}{\mathrm{SL}}
\newcommand{\Sp}{\mathrm{Sp}}
\newcommand{\Ep}{\mathrm{Ep}}
\newcommand{\OX}{\mathcal{O}(X)}

\newcommand{\SlC}{\mathrm{SL}_n(\mathbb{C})}
\newcommand{\SpC}{\mathrm{Sp}_{2n}(\mathbb{C})}
\newcommand{\EpC}{\mathrm{Ep}_{2n}(\mathbb{C})}
\newcommand{\Eu}[1]{\begin{pmatrix}
I&#1\\0&I
\end{pmatrix}}
\newcommand{\El}[1]{\begin{pmatrix}
I&0\\#1&I
\end{pmatrix}}

\newtheorem{theorem:introduction}{Theorem}
\newtheorem{theorem}{Theorem}[section]
\newtheorem*{theorem*}{Theorem}

\newtheorem{prop}[theorem]{Proposition}
\newtheorem{defi}[theorem]{Definition}

\theoremstyle{remark}
\newtheorem{Rem}[theorem]{Remark}

\author{Gaofeng Huang, Frank Kutzschebauch and Josua Schott \thanks{All three authors were partially supported by SNF-Grant 200021-207335.}}

\title{Factorization of Holomorphic Matrices and Kazhdan's property (T) }

\begin{document}
\maketitle

\begin{abstract}
    In this article we deduce some algebraic properties for the group $\Sp_{2n} (\mathcal{O} (X))$ of holomorphic symplectic matrices on a Stein space $X$: holomorphic factorization, exponential factorization, and Kazhdan's property (T).  In holomorphic factorization we combine a recent result of the third author and K-theory tools to give explicit bounds for the case when $X$ is one-dimensional or two-dimensional. Next we use them to find bounds for exponential factorization. As a further application, we show that the elementary symplectic group $\mathrm{Ep}_{2n}(\OX)$ admits Kazhdan's property (T). 
\end{abstract}

\tableofcontents
\section{Introduction}
In this paper we study algebraic properties of the general linear groups $\mbox{GL}_n (\mathcal{O} (X))$, the special linear groups $\mbox{SL}_n (\mathcal{O} (X))$ and the symplectic groups $\Sp_{2n} (\mathcal{O} (X))$ over the ring of holomorphic functions $\mathcal{O} (X)$ on a reduced Stein space $X$. Our first result is concerned with the classical K-theoretic question about the $\rm{K}_1$ groups of the ring $R=\mathcal{O} (X)$ and  bounded generation of the corresponding elementary subgroups with concrete bounds. Since we assume that the Stein space $X$ has finite dimension $n$ (defined as the complex dimension of the manifold $X \setminus \mathrm{Sing} (X)$, called the smooth part of $X$), these rings are interesting rings since they have finite Bass stable rank $\mathrm{bsr}( \mathcal{O} (X) ) = \lfloor \tfrac{1}{2}\dim X \rfloor +1 $ as established by Alexander Brudnyi in \cite{zbMATH07036259}. 
\begin{theorem:introduction}(Holomorphic Factorization) \label{theorem1}
    Let $X$ be a reduced Stein space of dimension $d$ and let $f: X \to \SpC$ be a null-homotopic holomorphic mapping with $n \ge 2$. Then
    \begin{enumerate}[(a)]
        \item There exists an upper bound $t=t(n,d)$ for the number of unitriangular factors in the factorization of $f$, considered as a holomorphic symplectic matrix in $\Sp_{2n}(\mathcal{O}(X))$. Moreover, $t(n,d)$ is smaller than or equal to the corresponding bound for $\Sl_{2}(\OX)$. 
        \item $t(n,1) = 4$ and $t(n,2) \le t(1,2)= 5$ for all $n \ge 2$.
    \end{enumerate}
\end{theorem:introduction}

The straightforward generalizations of our results to non-reduced Stein spaces can be achieved as in the aforementioned paper by Brudnyi. Here we consider only reduced Stein spaces. Our results rely on a combination of K-theoretic methods (Theorem \ref{theorem:tavgen}) and deep results from complex analysis, the solution to the so-called Gromov-Vaserstein problem (see Theorem \ref{theorem:IKLS}), about unitriangular (either upper triangular with 1 along the diagonal or lower triangular with 1 along the diagonal) factorization of holomorphic matrices. 
Since Theorem \ref{theorem:tavgen} is formulated in \cite{vavilov2011} only for the so-called UL factorization, i.e., for one with even number of factors, we give a complete proof for any (also odd) number of factors in section \ref{reduction}.

Next we consider our first application: the product of exponentials for the above mentioned groups $\mathrm{GL}_n (\mathcal{O} (X))$,  $\mathrm{SL}_n (\mathcal{O} (X))$,  $\mathrm{Sp}_{2n} (\mathcal{O} (X))$. We give a simple lower bound for all of these cases, prove the existence of an upper bound for the symplectic case, and give new upper bounds for the symplectic case for Stein spaces of dimension $1$ and dimension $2$. 
\begin{theorem:introduction}(Product of Exponentials) \label{theorem2}
    Let $X$ be a reduced Stein space of dimension $d$ and let $f:X\to \SpC$ be a null-homotopic holomorphic mapping with $n \ge 2$. Then
    \begin{enumerate}[(a)]
        \item There exist a natural number $e = e_{\mathrm{Sp}}(2n,d)$ and holomorphic mappings $A_1,...,A_e:X\to \mathfrak{sp}_{2n}(\mathbb{C})$ such that
        \begin{align*}
            f(x) = \exp(A_1(x))\cdots \exp(A_e(x)).
        \end{align*} 
        \item $ 2\le e_{\mathrm{Sp}}(2n,1)\le  e_{\mathrm{Sp}}(2n,2) \le 3 $.
    \end{enumerate}
\end{theorem:introduction}

We give another result in this direction, namely (Proposition \ref{proposition:eGL}) that the number of exponentials for the general linear group over the ring $\OX$ of holomorphic functions on a Stein space $X$  is at least $2$. For $2\times 2$ matrices this has been proved before by Studer and the second author in \parencite{zbMATH07202633}.

Our second application of the existence of uniform factorization is the fact that the path-connected component of the group $\Sp_{2n}(\mathcal{O} (X))$, which is equal to the elementary symplectic group $\mathrm{Ep}_{2n}(\OX)$, admits Kazhdan's
property (T) for $n\ge 2$. The corresponding result for the path-connected components of the groups
  $\Sl_n (\mathcal{O} (X))$, $n \geq 3$ is due to  Bj\"orn Ivarsson and the second author \cite{zbMATH06291038}. 
  \begin{theorem:introduction}(Kazhdan's Property (T)) \label{theorem3}
      Let $n \geq 2$ and let $X$ be a Stein manifold with finitely many connected components. Then $\mathrm{Ep}_{2n} (\OX) = (\mathrm{Sp}_{2n}(\OX))_0$ has Kazhdan's property (T).
  \end{theorem:introduction}
  
Since the solution to the Gromov-Vaserstein problem involves the so-called Oka principle, called the most beautiful principle in analysis by Ren\'e Thom, it is natural to compare the K-theoretic questions for the ring $\mathcal{O} (X) $ with the corresponding questions for the ring $\mathcal{C} (X)$ of continuous complex-valued functions on the Stein space $X$. More precisely, the Oka principle vaguely stated says that under certain conditions, the existence of a continuous solution implies the existence of a holomorphic solution. Let $t(n,d,\mathcal{C}, \mathcal{O})$ (see \cite{zbMATH06024090}) be the minimal number such that all null-homotopic holomorphic mappings, from a Stein space of dimension $d$ into $\mbox{SL}_n(\mathbb{C})$, factorize as a product of $t(n,d,\mathcal{C}, \mathcal{O})$ continuous unitriangular matrices
(starting with a lower triangular one) and let $t(n,d,\mathcal{O})$ be the minimal number that all null-homotopic holomorphic mappings, from Stein spaces of dimension
$d$ into $\mbox{SL}_n(\mathbb{C})$, factorize as a product of $t(n,d,\mathcal{O})$ holomorphic unitriangular matrices
(starting with a lower triangular one). In \cite{zbMATH06024090} it is proved 
that $t(2,1,\mathcal{C}, \mathcal{O})= t(2,1, \mathcal{O})=4$ and that if a fixed holomorphic matrix $A \in \mbox{SL}_2 (\mathcal{O} (X))$ factorizes
as a product of $N$ continuous unitriangular matrices, then it factorizes
as a product of $N+2$ holomorphic unitriangular matrices. Moreover, it is proved that the famous Cohn example \cite{zbMATH03232805}  $A_0 \in \Sl_2 (\mathcal{O} (\C^2))$ factorizes as a product of $4$ continuous 
 unitriangular matrices, but not less than $5$ holomorphic unitriangular matrices. Thus the question remained whether $t(2,2,\mathcal{C}, \mathcal{O})$ is equal to $4$ or to $5$. In the last section of this paper, we give an answer by proving
$t(2,2,\mathcal{C}, \mathcal{O})= t(2,2, \mathcal{O})=5$.

For matrices of bigger size it follows from \parencite{zbMATH06014072}, $\rm{bsr}(\OX) = 1$ when $\dim X = 1$ and \parencite{vavilov2011} that $t(n,1,\mathcal{C}, \mathcal{O})= t(n,1, \mathcal{O})=4$ for all $n \ge 3$. Note that $4 \le t(n,2,\mathcal{C}, \mathcal{O})\le  t(n,2, \mathcal{O})\leq 5$ (see \parencite[Remark 1.2]{zbMATH07036259}). Also if we denote the corresponding numbers for the symplectic group $\rm{Sp}_{2n}$ with the subscript symp, one has by Theorem \ref{theorem:tavgen}
\begin{align*}
    4 \le t_{\text{symp}}(2n,2,\mathcal{C}, \mathcal{O})\le  t_{\text{symp}}(2n,2, \mathcal{O}) \le 5
\end{align*} 
for all $n \ge 2$.

\section{Preliminaries}

\subsection{Elementary Generators}

Consider $n\times n$ matrices with entries in a commutative ring $R$ with $1$. Let $E_{ij}, i,j = 1, \dots, n$ be the matrix with $0$ everywhere except at the $(i,j)$ entry with $1$. Let $I$ be the $n \times n$ identity matrix. Then $I + r E_{ij}$ is upper diagonal with $1$ along the diagonal for $i < j, r \in R$, and lower diagonal with $1$ along the diagonal for $i > j, r \in R$. A product of matrices $I + r E_{ij}$ for $i < j, r \in R$ is upper triangular with $1$ along the diagonal, similarly for $i>j$ lower triangular. These upper and lower triangular matrices with $1$ along the diagonal are called \textit{unitriangular}. Let $\rm{E}_n(R)$ be the group generated by $I + r E_{ij}, i \neq j, r \in R$, and let $\Sl_n (R)$ be the set of matrices with determinant $1$. Over the complex numbers, $\rm{E}_n (\C) = \SlC$. 

Next, the symplectic group $\SpC$ can be represented as isometries of $\C^{2n}$ with respect to a nondegenerate, skew-symmetric bilinear form.  A convenient choice for the Gramian matrix of this bilinear form is 
\begin{align}
    J = \begin{pmatrix}    0 & I \\ -I & 0  \end{pmatrix},
\end{align}
where $0$ denotes the $n \times n$ zero matrix. We index rows and columns by $1, 2, \dots, n$,\\ $-1, -2, \dots, -n$.
In the block notation
\begin{align*}
    M = \begin{pmatrix}  A & B \\ C & D    \end{pmatrix} \in \SpC,
\end{align*}
the symplectic condition $M J M^T = J$ gives rise to three simple types of $J$-symplectic matrices: 

\begin{itemize}
    \item (i): $\begin{pmatrix}  I & B \\ 0 & I    \end{pmatrix}$, upper triangular with symmetric  $B = B^T$.
    \item (ii): $\begin{pmatrix}  I & 0 \\ C & I    \end{pmatrix}$, lower triangular with symmetric $C = C^T$.
    \item (iii): $\begin{pmatrix}  A & 0 \\ 0 & D    \end{pmatrix}$, block diagonal with invertible $A \in \rm{GL}_n(\mathbb{C})$ and $D = (A^{-1})^T$.
\end{itemize}

If block $A$ in type (iii) is upper triangular, then $D = (A^{-1})^T$ is lower triangular. In fact, $A$ and $D$ are simultaneously upper or lower triangular in another basis. This new basis can be obtained from the old one by reversing the order of the last $n$ basis elements, giving a Gramian matrix
\begin{align}
    \tilde{J} = \begin{pmatrix}    0 & L \\ -L & 0  \end{pmatrix}, 
\end{align}
where $L$ is the $n \times n$ matrix with $1$ along the skew-diagonal. Notice that symplectic matrices of type (i) and (ii) remain upper or lower triangular with respect to $\tilde{J}$, respectively. 


Over the ring $R$, the elementary $J$-symplectic generators for the elementary symplectic group $\Ep_{2n} (R)$ are similar to those for the special linear group. For example, type (i) corresponds to $I_{2n} + r (E_{i, -j} + E_{j, -i}), i \neq j$ and $I_{2n} + r E_{i, -i}$, while type (iii) corresponds to $I_{2n} + r (E_{ij} - E_{-j, -i}), i \neq j$ \parencite[186]{carter}. The subgroup $U$ of upper unitriangular $J$-symplectic matrices are generated by 
\begin{align*}
     \begin{pmatrix}  I & r (E_{ij} + E_{ji}) \\ 0 & I    \end{pmatrix}, \quad 
     \begin{pmatrix}  I & r E_{ii} \\ 0 & I    \end{pmatrix}, \quad
     \begin{pmatrix}  I + rE_{ij} & 0 \\ 0 & I - r E_{ji}    \end{pmatrix}, \quad i < j,\, r \in R. 
\end{align*}
Notice that the last form has upper triangular counterpart as $\tilde{J}$-symplectic matrices. Similarly one finds the generators for the subgroup $U^-$ of lower unitriangular $J$-symplectic matrices. From now on, we shall abbreviate $J$-symplectic as symplectic. Over the complex numbers, $\EpC = \SpC$. 

\subsection{Factorization of Holomorphic Mappings into \texorpdfstring{$\SpC$}{S}}

Let $X$ be a finite-dimensional Stein space and let $\mathcal{O}(X)$ be the ring of holomorphic functions on $X$. Then $\mathrm{Sp}_{2n}(\mathcal{O}(X))$ denotes the symplectic group with entries in $\mathcal{O}(X)$. Interpreting this group as holomorphic mappings from $X$ to $\mathrm{Sp}_{2n}(\mathbb{C})$, we denote by $(\mathrm{Sp}_{2n}(\mathcal{O}(X)))_0$ the path-connected component containing the identity. From \parencite[Theorem 1.1]{2207.05389} we cite the following
\begin{theorem}[IKLS Theorem]\label{theorem:IKLS}
Let $n$, $d$ be natural numbers and let $X$ be a reduced Stein space of dimension $d$. Then $\Ep_{2n}(\OX)=(\Sp_{2n}(\OX))_0.$
Moreover, there is a natural number $K(n,d)$, depending only on $n$ and $d$, such that each null-homotopic matrix $M$, i.e. $M\in (\mathrm{Sp}_{2n}(\OX))_0 $, can be written as a product of no more than $K(n,d)$ symplectic matrices of type (i) and (ii). \end{theorem}
\begin{Rem} \label{d1d2nullhomotop}
    Observe that elementary matrices in $\Ep_{2n}(\mathcal{O}(X))$ of type (i) and (ii) are null-homotopic, since we can multiply the off-diagonal entries by $t$, that is,
    $$ \Eu{tB}\quad \text{and}\quad \El{tC}.$$
    Therefore null-homotopy of $M \in \mathrm{Sp}_{2n}(\mathcal{O}(X))$ is a necessary condition. In the case where $X$ is contractible or a Stein space of dimension $1$ or $2$, every matrix is null-homotopic, that is, $(\mathrm{Sp}_{2n}(\OX))_0 $ and $\mathrm{Sp}_{2n}(\OX)$ coincide.
\end{Rem}

\begin{Rem}
    When $\dim X = 1$ (or $\dim X = 2$), we will see in Theorem \ref{unitrifactord1d2} that $4$ (or $5$) factors are sufficient for the unitriangular factorization, that is, factoring into the unipotent subgroups $U^{\pm}$. However, the set of elementary factors in Theorem \ref{theorem:IKLS} are taken from type (i) and (ii), and only the existence of an upper bound depending on $d$ and $n$ is guaranteed. For the same factorization over a field, $5$ factors are optimal \parencite{5FactorField}. Thus we find it natural to ask the following question.  \\
    
    \noindent \textbf{Open Problem}: Let $X$ be a Stein space of dimension $d$. Is there a bound for the optimal number $K(n,d)$ of factorizing a null-homotopic holomorphic mapping from $X$ to $\mathrm{Sp}_{2n}(\mathbb{C})$ into factors of type (i) and of type (ii), such that $K(n,d)$ is independent of $n$? 
\end{Rem}

\subsection{Reduction to Smaller Matrices} \label{reduction}

To obtain factorization estimates for holomorphic matrices of arbitrary size, we will make use of the Tavgen reduction. This appears in the setting of elementary Chevalley groups. Let $\Phi$ be a reduced irreducible root system of rank $l \geq 2$ and let $R$ be a commutative ring with $1$. We choose an order on $\Phi$ and a system of fundamental roots $\Pi = \{ \alpha_1, \alpha_2, \dots, \alpha_l \}$. Each root $\alpha \in \Phi$ is an integral sum of the fundamental roots
\begin{align*}
    \alpha = \sum_{i=1}^{l} k_i(\alpha) \alpha_i,
\end{align*}
where the integer coefficients $k_i(\alpha)$ are either all $\geq 0$ or all $\leq 0$. For $r = 1,l$, we define the following subsets of $\Phi$
\begin{align*}
    \Delta_r = \{ \alpha \in \Phi: k_r(\alpha) = 0 \}, \quad \Sigma_r = \{ \alpha \in \Phi: k_r(\alpha) > 0 \}, \quad \Sigma_r^- = \{ \alpha \in \Phi: k_r(\alpha) < 0 \}.
\end{align*}
$\Delta_r$ is itself a root system of rank $l -1$. On the level of Dynkin diagram, we obtain $\Delta_r$ from $\Phi$ by taking away the first ($r=1$) or the last ($r=l$) fundamental root. The elementary Chevalley group $E(\Phi,R)$ of type $\Phi$ over $R$ is generated by root subgroups $X_{\alpha}, \alpha \in \Phi$
\begin{align*}
    E(\Phi, R) =  \{ x_{\alpha}(r) \mid \alpha \in \Phi, r \in R \}.
\end{align*}
The positive unipotent subgroup $U(\Phi,R)$ is generated by the root subgroups of positive roots
\begin{align*}
    U(\Phi, R) = \{ x_{\alpha}(r) \mid \alpha \in \Phi^+, r \in R \}.
\end{align*}
Similarly, $U^-(\Phi, R) = \{ x_{\alpha}(r) \mid \alpha \in \Phi^-, r \in R \}$. The following theorem was originally proved by Oleg Tavgen and adapted in \parencite{vavilov2011}, where the number of factors is even. For our estimates, we need the same result allowing odd number of factors. We remark that the shape of the starting factor, upper or lower, is also immaterial. 

\begin{theorem}[Tavgen-VSS] \label{theorem:tavgen}
	Let $\Phi$ be a reduced irreducible root system of rank $l \geq 2$ and let $R$ be a commutative ring with $1$. Suppose that for subsystems $\Delta = \Delta_1, \Delta_l$ of rank $l-1$ the elementary Chevalley group $E(\Delta, R)$ admits the unitriangular factorization with $L$ factors
	\begin{align*}	
	E(\Delta, R) = U^-(\Delta, R)U(\Delta, R) \cdots U^{\pm}(\Delta, R).	
	\end{align*}
	Then the elementary Chevalley group $E(\Phi, R)$ admits the unitriangular factorization with the same number of factors
	\begin{align*}	
	    E(\Phi, R) =  U^-(\Phi, R)U(\Phi, R) \cdots U^{\pm}(\Phi, R).	
	\end{align*}
\end{theorem}

\noindent \textbf{Proof} 
    We take \begin{align*}	Y = U^-(\Phi, R)U(\Phi, R) \cdots U^{\pm}(\Phi, R).		\end{align*}
	$Y$ is a nonempty subset of $E(\Phi, R)$, in particular it contains $1$. Since the group $E(\Phi, R)$ is generated by the following  root elements $ X = \{ x_{\alpha}(r) \mid \alpha \in \pm \Pi, r \in R \} \subset E(\Phi, R).$
	Notice that the generating set $X$ is symmetric, i.e. $X^{-1} = X$. We claim that $x_{\alpha}(r) Y \subset Y$ for $\alpha \in \pm \Pi$: 
	Since $l \geq 2$, $\alpha$ lies in at least one of the subsystems $\Delta_1, \Delta_l$. Suppose that $\alpha$ belongs to $\Delta = \Delta_r$, then we consider the Levi decomposition 
	\begin{align*}	U(\Phi, R) = U(\Delta, R) \ltimes E(\Sigma, R), \qquad U^-(\Phi, R) = U^-(\Delta, R) \ltimes E(\Sigma^-, R), 	\end{align*}
	where $\Sigma = \Sigma_r$ and $E (\Sigma, R) = \langle x_{\alpha}(r) \mid \alpha \in \Sigma, r \in R \rangle$.
	Since $U^{\pm}(\Delta, R)$ normalizes $E^{\pm}(\Sigma, R)$ \parencite[Theorem 8.5.2]{carter}, we can rewrite $Y$ as
	\begin{align*}
		Y &= U^-(\Phi, R)U(\Phi, R) \cdots U^{\pm}(\Phi, R) \\	
		&=  U^-(\Delta, R) E(\Sigma^-, R) \, U(\Delta, R) E(\Sigma, R) \, \cdots \, U^{\pm}(\Delta, R) E(\Sigma^{\pm}, R) \\	
		&= \left( U^-(\Delta, R) U(\Delta, R) \cdots U^{\pm}(\Delta, R) \right) 
		   E(\Sigma^-, R) E(\Sigma, R) \cdots E(\Sigma^{\pm}, R) \\
		&= E(\Delta, R) \, E(\Sigma^-, R) E(\Sigma, R) \cdots E(\Sigma^{\pm}, R), 
	\end{align*}
	where the last step follows from the assumption. 
	For $\alpha \in \Delta$, $x_{\alpha}(r)$ is an element in $E(\Delta, R)$, hence $x_{\alpha}(r) Y \subset Y$. 
	
\noindent	\textbf{Lemma} Let $G$ be a group and $Y \subset G$ be a nonempty subset. Given a symmetric, generating subset $X$ of $G$, if     $X Y \subset Y$, 
	then $Y$ is the group $G$. \\
	The lemma shows that $Y = E(\Phi, R)$.   \hfill $\square$

\section{Number of Unitriangular Factors}

Here we restate two results from \parencite{zbMATH06024090} and use them to show part (b) of Theorem \ref{theorem1}.

\begin{theorem*}
    Let $X$ be a one-dimensional Stein space and let $f:X \to \mathrm{SL}_2(\mathbb{C})$ be a holomorphic mapping. Then there exist holomorphic mappings $g_1, \dots, g_4:X \to \mathbb{C}$ such that 
    \begin{align}\label{d1-4factor}
        f(x) = 
        \begin{pmatrix}
            1 & 0\\ g_1(x) & 1
        \end{pmatrix}
        \begin{pmatrix}
            1 & g_2(x)\\ 0 & 1
        \end{pmatrix}
        \begin{pmatrix}
            1 & 0\\ g_3(x) & 1
        \end{pmatrix}
        \begin{pmatrix}
            1 & g_4(x)\\ 0 & 1
        \end{pmatrix}.
    \end{align}
\end{theorem*}

\begin{theorem*}
    Let $X$ be a two-dimensional Stein space and let $f:X \to \mathrm{SL}_2(\mathbb{C})$ be a holomorphic mapping. Then there exist holomorphic mappings $g_1, \dots, g_5:X \to \mathbb{C}$ such that 
    \begin{align*}
        f(x) = 
        \begin{pmatrix}
            1 & 0\\ g_1(x) & 1
        \end{pmatrix}
        \begin{pmatrix}
            1 & g_2(x)\\ 0 & 1
        \end{pmatrix}
        \begin{pmatrix}
            1 & 0\\ g_3(x) & 1
        \end{pmatrix}
        \begin{pmatrix}
            1 & g_4(x)\\ 0 & 1
        \end{pmatrix}
        \begin{pmatrix}
            1 & 0\\ g_5(x) & 1
        \end{pmatrix}.
    \end{align*}
\end{theorem*}

\begin{theorem}\label{unitrifactord1d2}
Let $X$ be a reduced Stein space of dimension $1$ or $2$, and let $f:X \to \SpC$ be a holomorphic mapping. There exist holomorphic mappings $g_1, g_2, \dots, g_t:X \to \mathbb{C}^{n(2n+1)}$ such that 
\begin{align*}
    f(x) = M_1(g_1(x)) M_2(g_2(x)) \cdots M_t(g_t(x)),
\end{align*}
where $t = 4$ for $\dim X=1$ and $t \le 5$ for $\dim X=2$.
\end{theorem}

\noindent Here $M_j$ is defined as 
\begin{align*}
    M_j(g_j(x)) = 
    \begin{pmatrix}
        1 & & 0 \\
          & \ddots & \\
        g_j(x) & & 1  
    \end{pmatrix}
\end{align*}
for odd $j$ and 
\begin{align*}
    M_j(g_j(x)) = 
    \begin{pmatrix}
        1 & & g_j(x) \\
          & \ddots & \\
        0 & & 1  
    \end{pmatrix}
\end{align*}
for even $j$.

\noindent \textbf{Proof of Theorem \ref{unitrifactord1d2}} Let $f \in \Sp_{2n}(\OX)$. Then $f \in (\Sp_{2n}(\OX))_0 = \Ep_{2n}(\OX)$ by Remark \ref{d1d2nullhomotop} and Theorem \ref{theorem:IKLS}. Recall that $\Ep_{2n}(\OX)$ is the elementary Chevalley group $\mathrm{E}(\mathrm{C}_n, \OX)$. The above two results on $\Sl_2(\C)$ provide the base case with the root system $\mathrm{A}_1$ with the bound $t$. Then Tavgen reduction (Theorem \ref{theorem:tavgen}) guarantees that $t$ is also an upper bound for the number of unitriangular factors for $f \in \mathrm{Ep}_{2n}(\OX)$.   \hfill $\square$

\section{Number of Exponential Factors}
Let $\mathcal{A}$ be a m-convex Fréchet algebra with $1$ (see \parencite{michael}). The exponential of a $n \times n$ matrix $A$ is given by the exponential series. Denote by $\rm{Exp}_n(\mathcal{A})$ the subgroup of $\Gl_n(\mathcal{A})$ generated by exponentials and by $e(n,\mathcal{A})$ the minimal number such that any matrix in $\rm{Exp}_n(\mathcal{A})$ factorizes as a product of $e(n, \mathcal{A})$ exponentials. Let $t(n, \mathcal{A})$ be the minimal number such that any element in the elementary Chevalley group $E(\Phi, \mathcal{A}) \subset \Gl_n(\mathcal{A})$ factorizes as a product of $t(n,\mathcal{A})$ unitriangular matrices. When no such number exists, set $t(n,\mathcal{A}) = \infty$.

Observe that for a nilpotent matrix $N$, 
\begin{align*}
    \log (I + N) = \sum_{k=1}^{\infty } \frac{(-1)^k}{k} N^k
\end{align*}
is a finite sum. Thus every unipotent matrix $A$ can be written as the exponential of $\log A$. Also under conjugation an exponential remains an exponential, since $$B A B^{-1} = \exp(B \cdot \log A \cdot B^{-1})$$ for any $n \times n$ invertible matrix $B$. To estimate the number of exponential factors, we use the following trick presented in \parencite[Lemma 2.1]{zbMATH07459399} successively: Suppose we have a product of $3$ unitriangular matrices $U_1 U_2 U_3$ with alternate forms, e.g. $U_1, U_3$ are upper unitriangular, and $U_2$ lower unitriangular. Then 
\begin{align*}
    U_1 U_2 U_3 = (U_1 U_2 U_1^{-1}) (U_1 U_3) 
\end{align*}
is a product of $2$ exponentials. Because $U_2$ conjugated by $U_1$ is an exponential, while the product $U_1 U_3$ of upper unitriangular matrices is again upper unitriangular, in particular an exponential. Therefore, from $3$ unitriangular factors we get $2$ exponentials. An additional factor $U_4$ increases our estimate by one, but upon that adding a $5$th factor does not enlarge the number of exponentials
\begin{align*}
    U_1 U_2 U_3 U_4 U_5 = (U_1 U_2 U_1^{-1}) (U_1 U_3 U_4 U_3^{-1} U_1^{-1}) ( U_1 U_3 U_5 ). 
\end{align*}
The last factor has the same form as $U_5$. This leads us to the following estimate for $e(n, \mathcal{A})$.
\begin{prop}\label{proposition:e-t-estimate}
    \begin{align*}
    e(n, \mathcal{A})  \leq \lfloor \frac{1}{2} t(n,\mathcal{A}) \rfloor + 1.
    \end{align*}
\end{prop}

\begin{defi}
    For $\mathcal{A} = \OX$, denote by $e_{\Sl}(n, \OX)$ the minimal number such that any matrix in $ \Sl_n(\OX) \bigcap \rm{Exp}_n(\OX) = (\Sl_n(\OX))_0$ factorizes as a product of $e_{\Sl}(n, \OX)$ exponentials, with holomorphic mappings $X \to \mathfrak{sl}_n (\mathbb{C})$ as exponents. In this notation, $e_{\mathrm{SL}}(n,d)$ is the lowest upper bound for the number of exponential factors $e_{\mathrm{SL}}(n, \OX)$ over all Stein spaces $X$ of dimension $d$. The existence of this lowest upper bound was proved in \parencite{zbMATH07138526}. Similarly, $e_{\mathrm{Sp}} (2n, \OX)$ is the minimal number such that any matrix in $\Sp_{2n}(\OX) \bigcap \rm{Exp}_{2n}(\OX) =  (\Sp_{2n}(\OX))_0$ factorizes as a product of $e_{\Sp}(2n, \OX)$ exponentials, with holomorphic mappings $X \to \mathfrak{sp}_{2n} (\mathbb{C})$ as exponents. In this notation, $e_{\mathrm{Sp}}(2n,d)$ is the lowest upper bound for the number of exponential factors $e_{\Sp}(2n, \OX)$ over all Stein spaces $X$ of dimension $d$.
\end{defi}

\noindent \textbf{Proof of Theorem \ref{theorem2}} Apply Theorem \ref{theorem1} and Proposition \ref{proposition:e-t-estimate} with $\mathcal{A}=\Sp_{2n}(\OX).$ \hfill $\square$

\begin{Rem}
    For $f \in \Gl_n(\OX)$ null-homotopic, composition with the determinant $\det \circ f: X \to \C^*$ is also null-homotopic. Thus there exists a holomorphic function $g: X \to \C$ such that $\exp \circ g = \det \circ f$. This implies that 
    \begin{align*}
        e(n, \OX) \le e_{\mathrm{SL}}(n,\OX),
    \end{align*}
    since $\exp(-\frac{g}{n}  I_n) f$ is in $\mathrm{SL}_n (\OX)$ and 
    \begin{align*}
        \exp \left( \frac{g}{n}  I_n \right) \exp(A) = \exp \left( \frac{g}{n} A \right) .
    \end{align*}
\end{Rem}

\begin{prop} \label{proposition:eGL}
    Let $X$ be a Stein space with $\dim X > 0$ and let $n \ge 2$. Then $e(n,\OX) \ge 2$.
\end{prop}

\noindent \textbf{Proof}    \,
The proof is essentially the same as in \parencite{zbMATH07459399}. Let $X' \subset X$ be an irreducible component with $\dim X' > 0$. Then there exist two distinct points $x_1, x_2 \in X'$, we choose a holomorphic function $h \in \OX$ and $x_1, x_2 \in X'$ such that $h(x_1)=0, h(x_2) = 2 \pi i$. Set $g = \exp h$. Let
    \begin{align*}
        T = \begin{pmatrix}
            g & 1 \\
            0 & 1
        \end{pmatrix}.
    \end{align*}
    The same argument in \parencite{zbMATH07459399} shows that there does not exist $S \in M_2 (\OX)$ with $S^2 = T$ and in particular, $T$ does not have a logarithm. So $e(2, \OX) \ge 2$. For $n > 2$, fix $M \in \C \setminus \{0, 1 \}$ and set
    \begin{align*}
        T_n = \begin{pmatrix}
            M I_{n-2} & 0 \\
            0         & T
        \end{pmatrix}.
    \end{align*}
    Suppose that $T_n$ had a logarithm, then there would exist
    \begin{align*}
        S_n = \begin{pmatrix}
            L_1 & L_2 \\
            L_3 & L_4
        \end{pmatrix}
        \in M_n (\OX)
    \end{align*}
    with the same block partition as $T_n$ and such that $S_n^2 = T_n$. Then we have $S_n T_n = T_n S_n$, which implies that 
    \begin{align*}
        L_2 (T - M I_2) = 0 \quad \text{and} \quad (T - M I_2) L_3 = 0.
    \end{align*}
    On $X' \setminus (\exp h)^{-1} (M)$, $T - M I_2$ is invertible, so $L_2 = L_3 = 0$. By the identity theorem, $L_2$ and $L_3$ vanish on $X'$. But this would imply that $L_4^2 = T$, a contradiction. Hence  $e(n, \OX) \ge 2$ for all $n \ge 2$.   \hfill $\square$










\section{Kazhdan's Property (T)}
Let $G$ be a topological group, $K \subset G$, $\epsilon > 0$, $H$ a Hilbert space, and $(\pi, H)$ a continuous unitary $G$-representation. A vector $v \in H$ is called $(K, \epsilon)$-\textit{invariant} if $\| \pi (g)v - v \| < \epsilon \| v \|$ for all $g \in K$. We say that $G$ has \textit{Kazhdan's property (T)}, if there exist compact $K \subset G$ and $\epsilon > 0$ such that every continuous unitary $G$-representation with a $(K, \epsilon)$-invariant vector contains a nonzero $G$-invariant vector. $(K, \epsilon)$ is called a \textit{Kazhdan pair} for $G$. 
The symplectic group $\Sp_{2n}(R)$ over commutative ring $R$ is called \emph{boundedly elementary generated}, if there exists an integer $\nu$ such that every element is a product of at most $\nu$ elementary symplectic matrices. The following theorem is a symplectic version of \parencite[Theorem 2.6]{zbMATH06291038}. 

\hspace{5mm}

\noindent \textbf{Proof of Theorem \ref{theorem3}} 
    According to Theorem \ref{theorem:IKLS} and the fact that each factor of type (i) or (ii) is a product of at most $n(n+1)/2$ elementary symplectic generators, $\mathrm{Ep}_{2n} (\OX)$ is boundedly elementary generated. A finite set of holomorphic functions which generates a dense subring of $\OX$ can be constructed as follows. First, embed $X$ into $\C^N$ using Remmert's embedding theorem. Next, the set $S = \{ z_1, z_2, \dots, z_N, \sqrt{2}, i \}$ generates a dense subring of $\C [z_1, z_2, \dots, z_N ]$. Then, $\C [z_1, z_2, \dots, z_N ] \vert_X$ is dense in $\OX$ by the Oka-Weil theorem. Hence, $S$ generates a dense subring of $\OX$. Now \parencite[Theorem 1.1]{zbMATH02118941} gives a Kazhdan pair for $\Ep_{2n} (\OX)$.       \hfill  $\square$
    
\vspace{5mm}

From its definition, Kazhdan's property (T) is preserved under closure. As in \parencite[Theorem 3.1]{zbMATH06291038}, we can relate Kazhdan's property (T) of $\rm{Sp}_{2n} (\OX)$ to that of its quotient over the closure of the elementary group, with identical proof. 
    
\begin{theorem}
    Let $n \ge 2$ and let $X$ be a Stein space with finite embedding dimension. Then $\rm{Sp}_{2n} (\OX)$ has Kazhdan's property (T) if and only if 
    \begin{align*}
        \rm{Sp}_{2n}(\OX) / \overline{\mathrm{Ep}_{2n}(\OX)}
    \end{align*}
    has Kazhdan's property (T).
\end{theorem}


\section{Continuous versus Holomorphic Factorization}

In the following we present an example showing that $t(2,2,\mathcal{C}, \mathcal{O})$ = 5. Let $X$ be a two-dimensional Stein space, and let $f: X \to \Sl_2(\C)$ be a holomorphic mapping. Then $f$ can be written as a product of $5$ unitriangular holomorphic matrices by Theorem \ref{unitrifactord1d2}. In \parencite{zbMATH06024090} one finds an example which factorizes into $5$ unitriangular holomorphic factors and $4$ unitriangular continuous factors. Here we aim to give another example, which factorizes into $5$ unitriangular holomorphic factors \emph{and} $5$ continuous factors. To this end, we first study what it means to have a factorization with $4$ factors. Denote $f$ by $\begin{pmatrix}    a & b\\ c & d    \end{pmatrix}$, where $a, b, c, d \in \OX$ and $ ad - bc = 1$. If there exist $g_1, g_2, g_3, g_4: X \to \C$ so that $f$ can be decomposed as in (\ref{d1-4factor}), then the factorization reads
\begin{align*}
    \begin{pmatrix}
        1 & 0\\ g_1 & 1
    \end{pmatrix}
    \begin{pmatrix}
        1 & g_2\\ 0 & 1
    \end{pmatrix}
    \begin{pmatrix}
        1 & 0\\ g_3 & 1
    \end{pmatrix}
    \begin{pmatrix}
        1 & g_4\\ 0 & 1
    \end{pmatrix}
    = \begin{pmatrix}    a & b\\ c & d    \end{pmatrix}.
\end{align*}
Bring the first and the fourth factor to the right side, and carry out the multiplications 
\begin{align*}
    \begin{pmatrix}    1 + g_2 g_3 & g_2\\ g_3 & 1    \end{pmatrix} = 
    \begin{pmatrix}
        a & b - a g_4\\ c - a g_1 & - g_4 (c - a g_1) + d - b g_1
    \end{pmatrix}.
\end{align*}
In case $a \neq 0$, the first three equations read
\begin{align*}
    a = & 1 + g_2 g_3,  \\
    g_4  = &  \frac{1}{a} (b - g_2), \\
    g_1  = &  \frac{1}{a} (c - g_3), 
\end{align*}
and the fourth equation follows from the other three. If moreover $a \neq 1$, any choice of $g_3:\{ x \in X \mid a(x) \notin \{0,1\} \} \to \C^*$ gives a factorization in this part of $X$. The fiber of the fibration $f^* \Phi_4$ (see \parencite{zbMATH06024090}) over $\{ x \in X \mid a(x) \notin \{0,1\} \}$ is $\C^{*}$, where
\begin{align*}
    \Phi_4: \C^4 \to \Sl_2(\C), \,
    (z_1, z_2, z_3, z_4) \mapsto 
    \begin{pmatrix}
        1 & 0\\ z_1 & 1
    \end{pmatrix}
    \begin{pmatrix}
        1 & z_2\\ 0 & 1
    \end{pmatrix}
    \begin{pmatrix}
        1 & 0\\ z_3 & 1
    \end{pmatrix}
    \begin{pmatrix}
        1 & z_4\\ 0 & 1
    \end{pmatrix}.
\end{align*}

When $a = 0$, then
\begin{align*}
    1 + g_2 g_3 = 0, \quad g_2 = b, \quad g_3 = c, \quad 1 = - c g_4 + d - b g_1.
\end{align*}
Notice that $g_2$ and $g_3$ are prescribed as $b$ and $c$, respectively, and the fiber of $f^* \Phi_4$ here is $\C$. For $a = 1$, the fiber is the cross of axis. 

Consider the following holomorphic mapping $f: \C^2 \to \Sl_2(\C)$
\begin{align*}
    f(z, w) = \begin{pmatrix}    (zw-1)(zw-2) & (zw-1)z + (zw-2)z^2 \\ h_1(z,w) & h_2 (z,w)   \end{pmatrix},
\end{align*}
where the functions in the second row are chosen such that $f(z,w)$ has determinant $1$. The existence of such polynomial functions follows from Hilbert's Nullstellensatz, or if one is looking for holomorphic functions from a standard application of Theorem B. For this observe that the functions in the first row have no common zeros.

Suppose that there were continuous $g_1, g_2, g_3, g_4: \C^2 \to \C$ such that $f$ factorizes as in (\ref{d1-4factor}). Then on $zw = 1$, $g_2(z,w) = -z^2$ and on $zw = 2$, $g_2(z,w) = z$. Denote by $\xi_1, \xi_2$ the roots of $(D-1)(D-2)=1$, and choose a continuous curve $\gamma: [0,1] \to \C \setminus \{\xi_1, \xi_2\}$ such that $\gamma(0) = 1$ and $\gamma(1) = 2$. Then $g_2$ induces a family of continuous self-maps of $\C^*$
\begin{align*}
    F: [0,1] \times \C^* \to \C^*, (t, \theta) \mapsto g_2(\theta, \frac{1}{\theta}\gamma(t)).
\end{align*}
connecting  between $F(0,\theta) = -\theta^2$ and $F(1,\theta)=\theta$. But since these two self-maps of $\C^*$ have different degrees, we find a contradiction. 

\printbibliography

\end{document}